\documentclass[runningheads]{llncs}
\usepackage[T1]{fontenc}
\usepackage{graphicx}
\usepackage[pdfauthor   = {Mohamed\ Barakat\ and\ Marc\ Talleux\ and\ Fabian\ Zickgraf},
            pdftitle    = {Implementing\ the\ biset\ category\ of\ finite\ groups},
            pdfsubject  = {},
            pdfkeywords = {biset\ categories\; category\ of\ presheaves\; coequalizer\ completion\; coproduct\ completion\; ur-algorithms\; Schreier-Sims\ orbit\ algorithm},
            bookmarks=true,
            bookmarksopen=true,
            pagebackref=true,
            hyperindex=true,
            colorlinks=true,
            linkcolor=blue,
            citecolor=blue,
            filecolor=blue,
            urlcolor=blue,
            ]{hyperref}

\usepackage{amsmath}
\usepackage{amssymb}
\usepackage{mathtools}
\usepackage{mathrsfs}
\usepackage{latexsym}
\usepackage{xcolor}
\usepackage{tikz}
\usetikzlibrary{arrows,automata,positioning,backgrounds,calc,arrows.meta}
\usepackage{tikz-cd}
\usepackage{enumitem}
\usepackage{cases}
\usepackage{multicol}
\usepackage{float} 
\usepackage[normalem]{ulem} 
\usepackage[sort&compress,capitalise]{cleveref}
\usepackage{comment}
\usepackage{xspace}

\DeclareFontFamily{U}{dmjhira}{}
\DeclareFontShape{U}{dmjhira}{m}{n}{ <-> dmjhira }{}

\DeclareRobustCommand{\Y}{\text{\usefont{U}{dmjhira}{m}{n}\symbol{"48}}}

\newcommand{\CAP}{\textsc{Cap}\xspace}
\newcommand{\CompilerForCAP}{\ensuremath{\mathtt{CompilerForCAP}}\xspace}
\newcommand{\SkeletalFinSets}{\ensuremath{\mathtt{SkeletalFinSets}}\xspace}
\newcommand{\FinGSetsForCAP}{\ensuremath{\mathtt{FinGSetsForCAP}}\xspace}

\DeclareMathOperator\Nor{N}

\newcommand\N{\mathbb{N}}

\newcommand\C{\mathbf{C}}
\newcommand\D{\mathbf{D}}
\newcommand\op{\mathrm{op}}
\newcommand\E{\mathbf{E}}
\newcommand\T{\mathbf{T}}
\newcommand\V{\mathbf{V}}
\DeclareMathOperator\Hom{Hom}

\newcommand\LL{\mathscr{L}}
\newcommand\RR{\mathscr{R}}
\newcommand\UU{\mathscr{U}}

\newcommand{\rel}{\rightsquigarrow}

\DeclareMathOperator\id{id}
\DeclareMathOperator\Id{Id}

\DeclareMathOperator\coeq{coeq}
\newcommand\Coeq{\mathtt{Coeq}}

\newcommand\Set{\mathbf{Set}}
\newcommand\FinSet{\mathbf{FinSet}}
\newcommand\Cat{\mathbf{Cat}}
\newcommand\FinCat{\mathbf{FinCat}}
\newcommand{\FinBisetCat}{\mathbf{FinBisetCat}}
\newcommand{\FinBisetGrp}{\mathbf{FinBisetGrp}}
\newcommand{\FinBiset}{\mathbf{FinBiset}}
\newcommand{\Biset}{\mathbf{Biset}}
\newcommand\Ind{\mathbf{Ind}}
\newcommand\Inf{\mathbf{Inf}}
\newcommand\Def{\mathbf{Def}}
\newcommand\Res{\mathbf{Res}}
\newcommand\Iso{\mathbf{Iso}}
\DeclareMathOperator\stab{\mathrm{Stab}}
\DeclareMathOperator\PSh{\mathrm{PSh}}
\DeclareMathOperator\Func{\mathrm{Func}}
\DeclareMathOperator{\Aut}{Aut}

\newcommand{\CoequalizerOfIdentityAndAutomorphisms}{\mathtt{CoequalizerOfIdentity}\-\mathtt{AndAutomorphisms}}

\newcommand{\GFinSet}{\Gamma\textnormal{-}\mathbf{FinSet}}
\newcommand{\HFinSet}{H\textnormal{-}\mathbf{FinSet}}
\newcommand{\KFinSet}{K\textnormal{-}\mathbf{FinSet}}

\makeatletter
\def\slashedarrowfill@#1#2#3#4#5{%
  $\m@th\thickmuskip0mu\medmuskip\thickmuskip\thinmuskip\thickmuskip
   \relax#5#1\mkern-7mu%
   \cleaders\hbox{$#5\mkern-2mu#2\mkern-2mu$}\hfill
   \mathclap{#3}\mathclap{#2}%
   \cleaders\hbox{$#5\mkern-2mu#2\mkern-2mu$}\hfill
   \mkern-7mu#4$%
}
\def\rightslashedarrowfill@{%
  \slashedarrowfill@\relbar\relbar\mapstochar\rightarrow}
\newcommand\xslashedrightarrow[2][]{%
  \ext@arrow 0055{\rightslashedarrowfill@}{#1}{#2}}
\makeatother

\newcounter{CntWP}                                      
\newcounter{CntT}[CntWP]                        

\usepackage[draft]{optional}
\usepackage[colorinlistoftodos,prependcaption,textsize=tiny]{todonotes}

\newcommand{\plan}[1]{}
\newcommand{\BA}[1]{}
\newcommand{\MB}[1]{}
\opt{draft}{
   \renewcommand{\plan}[1]{{\color{blue}{#1}}\PackageWarning{TODO}{Plan: #1}}
   \renewcommand{\BA}[1]{\todo[color=orange!30]{BA: #1} \PackageWarning{TODO}{BA: #1}}
   \renewcommand{\MB}[1]{\todo[color=green!30]{MB: #1}\PackageWarning{TODO}{MB: #1}}
}

\begin{document}
\title{Implementing the biset category of finite groups}
%
%
\author{Mohamed Barakat\inst{1}
\and
Marc Talleux\inst{2}
\and
Fabian Zickgraf\inst{1}
}
%
%
\institute{Department of mathematics, University of Siegen, 57068 Siegen, Germany\\
\email{mohamed.barakat@uni-siegen.de}, \email{fabian.zickgraf@uni-siegen.de}\\
\and
LAMFA, Université de Picardie Jules Verne, 80039 Amiens, France\\
\email{marc.talleux@u-picardie.fr}}
\maketitle              

{\centering\footnotesize \em Dedicated to Irene Schreier-Scott and Dana Scott\par}

\begin{abstract}
We describe an implementation of the biset category of finite groups as a tower of standard categorical constructions, all of which are implemented in the software project \CAP for algorithmic category theory.
In particular, we describe the composition of bisets as a composition in a Kleisli category of some biadjunction monad.
This composition relies on the universal property of the coequalizer completion of a group viewed as a groupoid on one object.
Expressing this universal property offers an elegant categorical interpretation of the Schreier--Sims orbit algorithm.
Indeed, the implementation relies on every aspect of the algorithm.

\end{abstract}

\section{Introduction} \label{sec:classical_bisets}

Throughout this paper $G,H,K$ will denote finite groups.
Many of the constructions do not need the finiteness assumption.\footnote{The finiteness assumption makes the sum in the Mackey formula \eqref{eq:Mackey} finite.}

An $(H,G)$-biset is a set $X$ equipped with a left $H$-action and a commuting right $G$-action:
\begin{equation} \label{hxg}
  h (x g) = (h x) g.
\end{equation}

The \emph{biset category $\FinBiset$ of finite groups} \cite{bouc2010biset} has finite groups as objects and isomorphism classes of $(H,G)$-bisets as morphisms in $\FinBiset(G,H)$.
Given two bisets ${}_H X_G$ and ${}_K Y_H$ their composition $({}_H X_G, {}_K Y_H) \mapsto {}_K Y_H \circ {}_H X_G$ is defined by
\[
  {}_K Y_H \circ {}_H X_G \coloneqq {}_K Y \times_H X_G \mbox{,}
\]
where $Y\times_H X$ is the set of orbits of the diagonal action of $H$ on $Y \times X$ given by $h (y, x) \coloneqq (y h^{-1}, h x)$.
In \Cref{sec:goursat_bouc} we introduce the notion of a transitive biset and briefly review the Goursat-Bouc decomposition which decomposes any transitive biset as the composition of elementary transitive bisets.

The goal of this paper is to describe our \CAP implementation of the biset category $\FinBiset$ of finite groups.
\CAP is an open source software project for algorithmic category theory \cite{GPSSyntax,GP_Rundbrief,PosCCT+arXiv}, which allows the construction of categories as towers of standard categorical constructions.
It turns out that $\FinBiset$ is a full subcategory of the Kleisli category of the biadjunction monad of the finite colimit completion of small categories, where this completion is viewed as a $2$-functor, left biadjoint to the functor forgetting colimits.
The finite colimit completion of a category is the coequalizer completion of the finite coproduct completion.
We explain this in \Cref{sec:tower}.

In the case of a group $G$, viewed as a groupoid on one object, these two completions can be switched due to the orbit theorem.
The coequalizer completion of $G$ is nothing but the category of transitive left $G$-sets.
Expressing the universal property of the coequalizer completion amounts to a nontrivial application of the Schreier--Sims orbit algorithm.
This turns out to be the crucial step in describing the composition of bisets as the composition in the above-mentioned Kleisli category.
This is the content of \Cref{sec:implementation}.

\section{The Goursat-Bouc decomposition} \label{sec:goursat_bouc}

Let ${}_H X_G$ be an $(H,G)$-biset.
Each $g \in G$ defines, due to \eqref{hxg}, an $H$-equivariant automorphism $\chi_g\colon X \to X, x \mapsto x g$ of ${}_H X$, inducing a group homomorphism
\begin{equation} \label{eq:chi}
  \chi\colon G \to \Aut({}_H X), g \mapsto \chi_g \mbox{,}
\end{equation}
where the composition in $\Aut({}_H X)$ is the precomposition $\chi_g \chi_{g'} = \chi_{g g'}$.
In particular, the biset $X$ gives rise to the following right action of $G$ on the set $H \!\setminus\! X$ of $H$-orbits:
\begin{equation} \label{Hxg}
  H \!\setminus\! X \times G \to H \!\setminus\! X, (H \cdot x, g) \mapsto H \cdot x g \mbox{.}
\end{equation}
The biset ${}_H X_G$ is called \emph{transitive} if this action is transitive.

Such a biset can equally be viewed as a left $(H \times G)$-set, where the action is defined by $(h, g) x \coloneqq h x g^{-1}$.
The biset ${}_H X_G$ is transitive iff this action is transitive.
In fact, the transitive $(H,G)$-bisets are in bijection with the conjugacy classes of subgroups of $H\times G$:
A subgroup $D \leq H \times G$ gives rise to a transitive $(H,G)$-biset $(H \times G) / D$ under the action $h \cdot (x, y) D \cdot g \coloneqq (hx,g^{-1}y)D$.
One can recover $D$ as the stabilizer $D = \stab_{H \times G}( (1,1)D )$ of $(1,1)D \in (H \times G) / D$.
Conversely, for a transitive $(H,G)$-biset ${}_H X_G$ the stabilizer
\[
  D \coloneqq \stab_{H \times G}( x ) \coloneqq \{ (h, g) \in H\times G \mid h x = x g \}
\]
of an $x \in X$ is a subgroup of $H \times G$ and ${}_H X_G \cong (H \times G)/D$ as $(H,G)$-bisets.

Let $\Gamma$ be a group, $U \leq \Gamma$ a subgroup of $\Gamma$, $N \unlhd \Gamma$ a normal subgroup, and $\varphi\colon \Gamma' \to \Gamma$ an isomorphism of groups.
Denoting the factor group by $Q \coloneqq \Gamma/N$, we define the five \emph{elementary} transitive bisets
\begin{align*}
  \Ind_U^\Gamma \coloneqq {}_\Gamma \Gamma_U,
  \Res_U^\Gamma \coloneqq {}_U \Gamma_\Gamma,\,\,\,\,
  \Inf_Q^\Gamma \coloneqq {}_\Gamma Q_Q,
  \Def^\Gamma_Q \coloneqq {}_Q Q_\Gamma,\,\,\,\,
  \Iso(\varphi) \coloneqq {}_\Gamma \Gamma_{\Gamma'} \mbox{,}
\end{align*}
where maybe the only nonobvious action is the right action of $\Gamma'$ on $\Gamma$ via $\varphi$.
It follows that any transitive $(H,G)$-biset ${}_H X_G \cong (H \times G)/D$ can be written as the Goursat-Bouc decomposition \cite[Section~2.3]{bouc2010biset}:
\begin{align} \label{eq:fivefold}
  {}_H X_G \cong (H\times G)/D \cong \Ind^H_{P_1} \circ \Inf^{P_1}_{P_1/K_1} \circ \Iso(\varphi) \circ \Def^{P_2}_{P_2/K_2} \circ \Res^G_{P_2} \mbox{,}
\end{align}
where $P_1 = P_1(D)$, $P_2 = P_2(D)$ denote the projections of $D$ onto $H$ and $G$, respectively,
\begin{align*}
  K_1 &= K_1(D) \coloneqq \{h \in H \mid (h,1) \in D \} \unlhd P_1, \\
  K_2 &= K_2(D) \coloneqq \{g \in G \mid (1,g) \in D \} \unlhd P_2
\end{align*}
are subnormal subgroups, and $\varphi\colon P_2/K_2 \xrightarrow{\sim} P_1/K_1$ the isomorphism defined by Goursat's Lemma.
One can recover $D$ (and hence the transitive biset) from the two subquotients and the isomorphism via the equality
\begin{equation} \label{eq:D}
  D = \{\ (h, g) \in P_1 \times P_2 \mid h K_1 = \varphi( g K_2 )\ \} \mbox{.}
\end{equation}

The only nontrivial relation among elementary bisets is given by the \emph{Mackey formula}, expressing the $(V,U)$-biset $\Res^\Gamma_V \circ \Ind^\Gamma_U$ for two subgroups $U,V \leq \Gamma$ as the sum (=coproduct) of transitive $(V,U)$-bisets
\begin{align} \label{eq:Mackey}
  \Res^\Gamma_V \circ \Ind^\Gamma_U
  \cong
  \sum_{x \in [V \setminus \Gamma / U]} \Ind^V_{V \cap {}^x U} \circ \Iso(\gamma_x) \circ \Res^U_{V^x \cap U} \mbox{,}
\end{align}
with isomorphisms $\gamma_x\colon V^x \cap U \xrightarrow{\sim} V \cap {}^x U$ for each representative $x$ of a double coset in $V \setminus\! \Gamma / U$.

\section{Bisets from a categorical point of view} \label{sec:tower}

\subsection{Bisets for finite categories} \label{sec:Biset}

Let $\FinCat$ denote the $2$-category of finite categories, naturally enriched over the elementary topos $\FinSet$ of finite sets.
Let $\C, \D, \E \in \FinCat$ denote finite categories.

Define the category $\FinBisetCat$ with the same objects as $\FinCat$ but with functors $\C \to \D$ replaced by profunctors $X\colon \C \xslashedrightarrow{} \D$, i.e., functors $X\colon \D^\op \times \C \to \FinSet$.
For the wide range of applications of bisets for categories see \cite{web2023}.
The composition of profunctors is usually described by a coend-construction, which we do not know how to make algorithmically efficient.
Alternatively, one can use the $2$-adjunction in $\Cat$ to rewrite $X$ as a functor
\[
  X\colon \C \to \PSh(\D) \mbox{,}
\]
where $\PSh(\D) \coloneqq \PSh(\D, \FinSet) \coloneqq \Func(\D^\op, \FinSet)$ denotes the elementary topos of finite-set-valued presheaves  on $\D$.
Now in order to compose two functors $X\colon \C \to \PSh(\D)$ and $Y\colon \D \to \PSh(\E)$ to a functor $XY\colon \C \to \PSh(\E)$ we first need to lift $Y$ to a (finitely cocontinuous) functor $\widehat{Y}\colon \PSh(\D) \to \PSh(\E)$.
We call this composition $XY \coloneqq X \widehat{Y}$ the \emph{Kleisli composition}\footnote{The similarity with the composition in the Kleisli category is not a coincidence. Since $\PSh(-,\Set)$ is a model for the colimit completion of small categories, $\PSh(-,\Set)$ is part of a biadjunction and $\Biset$ is the Kleisli category of the adjunction monad.}.
The first obvious advantage of this approach is that the $2$-category of presheaf categories and cocontinuous functors $\widehat{Y}\colon \PSh(\D) \to \PSh(\E)$ is a strictification of the weak $2$-category of categories and profunctors $Y\colon \D \xslashedrightarrow{} \E$.

\subsection{The universal property of $\PSh$}

A second advantage is the insight that the composition of two profunctors can be reduced to the ur-algorithm\footnote{By an \emph{ur-algorithm} we refer in \CAP to an algorithm implementing the counit (=evaluation) of a pair $\LL \dashv \RR$ of biadjoint $2$-functors, expressing the universal property of $\LL$.
} expressing the universal property of the left biadjoint $2$-functor $\LL = \PSh$:
Given a finitely cocomplete category $\T$ and a functor $Y\colon \D \to \T$, determine a finitely cocontinuous lift $\widehat{Y}\colon \PSh(\D) \to \T$.
Indeed, the so-called co-Yoneda lemma states that each presheaf $\Omega \in \PSh(\D)$ can be ``resolved'' as the colimit of a diagram of representable presheaves $\Hom_\D(-,c)$, while the fully faithful Yoneda functor $\Y\colon \D \hookrightarrow \PSh(\D)$ allows one to consider this diagram as a diagram in $\D$, but where its colimit is taken in $\PSh(\D)$ and coincides with $\Omega$.
By applying $Y\colon \D \to \T$ to the resolving diagram of the presheaf $\Omega$ in $\D$ we get a diagram in $\T$, the colimit of which is the desired $\widehat{Y}(\Omega)$.
The so constructed finitely cocontinuous functor $\widehat{Y}\colon \PSh(\D) \to \T$ is the extension of $Y$ along $\Y$, expressing the \emph{universal property} of $\PSh(\D)$ as a model for the finite colimit completion of $\D$.
In the rest of this subsection we will explain how one can realize the finite colimit completion as the coequalizer completion $\Coeq(\sqcup \D)$ of the finite coproduct completion $\sqcup \D$ establishing the equivalence
\begin{equation} \label{PSh_as_colimit_completion}
  \PSh(\D) \simeq \Coeq(\sqcup \D) \mbox{.}
\end{equation}

The algorithmic value of this presheaf-approach to profunctors thus strongly depends on the efficiency of the ur-algorithm realizing the universal property of $\PSh(\D)$, which it turn depends on the size of the resolving diagram of the presheaf $\Omega$.
A resolving diagram $Q= Q_\Omega$ of $\Omega$ is formally a quiver with vertices being objects in $\D$ and arrows decorated by morphisms in $\D$.
Such a quiver $Q$ can itself be expressed (in \CAP) as a presheaf $Q \in \PSh\left( V \overset{s}{\underset{t}{\rightrightarrows}} A, \sqcup \D\right)$, i.e., a coequalizer pair of morphisms $Q(V) \overset{Q(s)}{\underset{Q(t)}{\leftleftarrows}} Q(A)$ in $\sqcup \D$.
By the universal property of the finite coproduct completion $\sqcup \D$, the Yoneda embedding $\Y$ extends to a full embedding $\sqcup \D \hookrightarrow \PSh(\D)$ recovering $\Omega$ as the coequalizer in $\PSh(\D) \supset \sqcup \D$:
\begin{equation} \label{eq:colim}
  \Omega \overset{\pi_Q}{\twoheadleftarrow} Q(V) \overset{Q(s)}{\underset{Q(t)}{\leftleftarrows}} Q(A) \mbox{.}
\end{equation}

A canonical, yet impractical such diagram $Q$ is given by Grothendieck's category of elements.
Each object $d \in \D$ appears $|\Omega(d)|$ times as a vertex and each morphism $\varphi\colon d \to d'$ appears $|\Omega(d')|$ times as an arrow in the diagram, i.e.,
$ Q(V) = \sqcup_{d \in \D} |\Omega(d)| \otimes d,\, Q(A) = \sqcup_{\varphi\colon d \to d' \in \D} |\Omega(d')| \otimes d'$.
The description of $Q(s)$ and $Q(t)$ is straightforward.
Although the number of arrows in $Q(A)$ can be reduced by only considering generating morphisms of $\D$, the potentially huge number of vertices renders this description algorithmically useless.
Finding efficient resolving diagrams $Q_\Omega$ to efficiently realize the universal property of $\PSh$ and use it to compose bisets of finite categories are left to future work.

\subsection{The category $\GFinSet$ as a tower}

Let $\Gamma$ be a finite group, viewed as a groupoid on one object $*$.
The presheaf category $\PSh(\Gamma)$ of $\Gamma$ with values in finite sets is just another encoding of the category $\GFinSet$ of finite \emph{left} $\Gamma$-sets.
Hence, by \eqref{PSh_as_colimit_completion}, we get the chain of equivalences:
\[
  \GFinSet \simeq \PSh(\Gamma) \simeq \Coeq(\sqcup \Gamma) \mbox{.}
\]
An elegant categorical interpretation of the orbit theorem for (finite) groups is that both completions in the categorical tower can be switched (which already fails for the cyclic monoid $\langle e | e^2 = e \rangle$ of order $2$), i.e.,
\begin{align} \label{eq:tower}
  \GFinSet \simeq \Coeq(\sqcup \Gamma) \simeq \sqcup \Coeq(\Gamma) \mbox{,}
\end{align}
where the coequalizer completion $\Coeq(\Gamma)$ is (equivalent to) the category of transitive left $\Gamma$-sets.
This flipped categorical tower will be the basis of the \CAP implementation of $\GFinSet$ in \Cref{sec:Gamma-FinSet_as_a_tower_implementation}.

\subsection{The biset category of finite groups} \label{sec:Biset_for_groups}

From now on, we will denote the full subcategory of $\FinBisetCat$ on finite groups (viewed as one-object groupoids) by $\FinBisetGrp$ or simply by $\FinBiset$.
It is easy to see that this is equivalent to the definition from \Cref{sec:classical_bisets}, i.e.,
\[
  \FinBiset(G,H) \cong \{ {}_H X_G \mid {}_H X_G \mbox{ an } (H,G)\mbox{-biset} \} \mbox{,}
\]
where we identify an $(H,G)$-biset ${}_H X_G$ with a functor
\begin{equation} \label{eq:biset}
  X\colon G \to \HFinSet \simeq \PSh(H) \simeq \Coeq(\sqcup H) \simeq \sqcup \Coeq(H) \mbox{.}
\end{equation}
Note that the Kleisli composition $X Y \in \FinBiset(G,K)$ of $X \in \FinBiset(G,H)$ and $Y \in \FinBiset(H,K)$ now corresponds to the post-composition ${}_K Y_H \circ {}_H X_G$.
The Goursat-Bouc decomposition \eqref{eq:fivefold} shows that $\FinBiset$ is generated (as an $\N$-linear category) by the elementary transitive bisets, subject to various commutativity relations and the Mackey relation \eqref{eq:Mackey}.

\section{Implementation in \CAP} \label{sec:implementation}

As explained in \Cref{sec:Biset} and detailed in \Cref{fig:Kleisli} below, the composition $XY$ of bisets reduces to the Kleisli composition $X \widehat{Y}$, where $\widehat{Y}\colon \HFinSet \to \KFinSet$ is the finitely cocontinuous lift of $Y\colon H \to \KFinSet$ along the Yoneda embedding $H \hookrightarrow \HFinSet$.
In \Cref{sec:Gamma-FinSet_as_a_tower_implementation} we explain how to construct the categorical tower $\sqcup \Coeq(\Gamma) \simeq \GFinSet$ \eqref{eq:tower} and in \Cref{sec:SchreierSim} how to compute the intermediate lift $\widetilde{Y}\colon \Coeq(H) \to \sqcup \Coeq(K)$ (see \Cref{fig:Kleisli}), which is the algorithmically crucial part of the finitely cocontinuous lift $\widehat{Y}\colon \sqcup \Coeq(H) \to \sqcup \Coeq(K)$.

\begin{figure}
\begin{center}
\begin{tikzpicture}[scale=0.98]
  \node (HFinSet) {$\HFinSet \simeq$};
  \node (UCoeqH) at ($(HFinSet)+(2,0)$) {$\sqcup \Coeq(H)$};
  \node (CoeqH) at ($(UCoeqH)+(0,-1)$) {$\Coeq(H)$};
  \node (H) at ($(CoeqH)+(0,-1)$) {$H$};
  \node (KFinSet) at ($(UCoeqH)+(3.5,0)$) {$\KFinSet$};
  \node (UCoeqK) at ($(KFinSet)+(2,0)$) {$\simeq \sqcup \Coeq(K)$};
  \node (G) at ($(H)+(-5.5,0)$) {$G$};
  
  \draw[right hook->] (H) -- (CoeqH);
  \draw[right hook->] (CoeqH) -- (UCoeqH);
  \draw[->] (G) -- node[midway, fill=white]{$X$} (HFinSet);
  \draw[->] (H) -- node[midway, fill=white]{$Y$} (KFinSet);
  \draw[thick,->] (CoeqH) -- node[midway, fill=white]{$\widetilde{Y}$} (KFinSet);
  \draw[->] (UCoeqH) -- node[midway, fill=white]{$\widehat{Y}$} (KFinSet);
  \draw[->] (G) [rounded corners=10pt] -- ($(HFinSet)+(-1,0.7)$) -- ($(UCoeqH)+(1,0.7)$) node[midway, fill=white]{$XY \coloneqq X \widehat{Y}$} -- (KFinSet);
\end{tikzpicture}
\end{center}
\caption{The Kleisli composition $X \widehat{Y}$ computes the composition $XY$ of bisets}
\label[figure]{fig:Kleisli}
\end{figure}

\subsection{Implementing the category $\GFinSet$ as a tower} \label{sec:Gamma-FinSet_as_a_tower_implementation}

\subsubsection{The coequalizer completion of $\Gamma$.}

A \emph{skeletal} model of the category $\Coeq(\Gamma)$ of transitive \emph{left} $\Gamma$-sets can be easily constructed:
Let $\{U_c\}_{c=1}^\ell$ be the set of representatives of the $\ell$ many conjugacy classes of subgroups of $\Gamma$ with $c \leq d \implies |U_c| \leq |U_d|$, in particular $U_1 = \{1_\Gamma\}$ (the trivial subgroup) and $U_\ell = \Gamma$.
Define the skeletal category $\Coeq(\Gamma)$ with objects $c \in \{1, \ldots, \ell \}$, corresponding to the \emph{transitive} $\Gamma$-sets $\{ \Gamma / U_c \}_{c=1}^\ell$.
In particular, $1$ now stands for the principal $\Gamma$-set.
A morphism $\gamma\colon s \to t$ is given by a group element $\gamma \in \Gamma$ satisfying
\begin{equation} \label{eq:mor_in_Coeq}
  \gamma^{-1} U_s \gamma \eqqcolon U_s^\gamma \leq U_t \mbox{,}
\end{equation}
where two parallel morphisms $\gamma, \gamma'\colon s \to t$ are considered equal iff $\gamma^{-1} \gamma' \in U_t$.
The identity on $c$ is $1_\Gamma\colon c \to c$ and the composition $s \xrightarrow{\gamma} c \xrightarrow{\gamma'} t$ is given by the group multiplication $s \xrightarrow{\gamma \gamma'} t$.
In $\Coeq(\Gamma)$ all morphisms are epis, and the classes of monos, endomorphisms, and automorphisms coincide.
The embedding $\Gamma \hookrightarrow \Coeq(\Gamma)$ sends the unique object $*$ to $1$ (corresponding to the principal $\Gamma$-set $\Gamma / U_1$) and each $\gamma \in \Gamma$ to the automorphism $\gamma\colon 1 \to 1$.

Indeed, the category $\C = \Coeq(\Gamma)$ has all coequalizers.
However, in what follows we will only need specific coequalizers, namely the coequalizer $\coeq_c(\alpha) = \coeq(\id_c, \alpha_1, \ldots, \alpha_k)$ of the identity $\id_c$ of an object $c \in \{1, \ldots, \ell \}$ and $k$ many automorphisms $\alpha = (\alpha_1, \ldots, \alpha_k)$  of $c$, where $\alpha_r \colon c \to c$ and $k \in \N$.
To compute $\coeq_c(\alpha)$, first define the subgroup $V \coloneqq \langle U_c, \alpha_1, \ldots, \alpha_k \rangle \leq \Gamma$.
Then compute $\gamma \in \Gamma$ such that $V^\gamma = U_d$ for some $d \in \{1, \ldots, \ell\}$.
It follows that, $\coeq_c(\alpha) = d$ and, using \eqref{eq:mor_in_Coeq} and the definition of $V$, $\pi = \gamma\colon c \to d$ is the projection onto the coequalizer.
In the very special case that $U_c$ is self-normalizing it follows that each $\alpha_r = \id_c$, hence $V = U_d = U_c$ and the projection $\pi$ is given by $1_\Gamma\colon c \to c$.
Which $U_c \leq \Gamma$ is self-normalizing can be read off the diagonal of the table of marks of $\Gamma$.
It remains to explain how to compute the colift along the coequalizer:
Given an epimorphism $\tau\colon c \to t$ with $\tau = \alpha_r \tau$ for all $r = 1, \ldots, k$, the colift along the projection $\pi = \gamma\colon c \to d$ is simply given by $\gamma^{-1}\tau\colon d \to t$.

The universal property of $\Coeq(\Gamma)$ as a coequalizer completion of the one-object groupoid $\Gamma$ is easy to describe:
Given a functor $Y\colon \Gamma \to \T$ into a category $\T$ with coequalizers, the lift $\widetilde{Y}\colon \Coeq(\Gamma) \to \T$ is constructed by the following simple ur-algorithm:
Write each object $c \in \Coeq(\Gamma)$ as the coequalizer of the identity $1_\Gamma\colon 1 \to 1$ (of the principal $\Gamma$-set) and the automorphisms $u_{c 1}, \ldots, u_{c k_c}\colon 1 \to 1$, where $\{ u_{c 1}, \ldots, u_{c k_c} \}$ is any generating set of $U_c$.
Now define $\widetilde{Y}(c)$ as the coequalizer in $\T$ of the identity $Y(1_\Gamma) = \id_{Y(*)}$ and the automorphisms $Y(u_{c 1}), \ldots, Y(u_{c k_c})\colon Y(*) \to Y(*)$.
Since such a specific coequalizer is central for lifting an action of a group $\Gamma$ to a functor $\Coeq(\Gamma) \to \T$ (and then to a functor $\GFinSet = \sqcup \Coeq(\Gamma) \to \T$), we give it a special name in \CAP:
\begin{equation} \label{eq:coequalizer_of_identity_and_automorphisms}
  \coeq_t(\alpha) \coloneqq \CoequalizerOfIdentityAndAutomorphisms( t, (\alpha_1, \ldots, \alpha_k) ) \mbox{.}
\end{equation}
In our case $t = Y(*) \in \T$ and $\alpha_r = Y(u_{c r}) \in \Aut_\T(t)$ for $r = 1, \ldots, k = k_c$.

\subsubsection{The finite coproduct completion of a skeletal object-finite category.}

Let $\C$ be a category with finitely many objects $\{C_1, \ldots, C_\ell\}$ (in our application $\C = \Coeq(\Gamma)$ with objects $\{1, \ldots, \ell\} \equiv \{\Gamma/U_1, \ldots, \Gamma/U_\ell\}$).
If $\C$ is even skeletal, then we can choose a model for the coproduct completion\footnote{The constructor $\mathtt{FiniteStrictCoproductCompletionOfObjectFiniteCategory}$ was implemented in the \CAP-based package $\mathtt{FiniteCocompletions}$ \cite{FiniteCocompletions}.}
\[
  \sqcup \C  \coloneqq \mathtt{FiniteStrictCoproductCompletionOfObjectFiniteCategory}( \C ) \mbox{,}
\]
which is both skeletal and strict:
An object $T \coloneqq \coprod_{c=1}^\ell C_c^{\sqcup m_c} \in \sqcup \C$ can be identified with the $\ell$-tuple $(m_1, \ldots, m_\ell)$ of multiplicities of the objects $C_c$ in the coproduct decomposition of $T$.
In this skeletal data structure the coproduct of two $\ell$-tuples is just their componentwise summation.

We will refrain from describing the data structure of a general morphism in $\sqcup \C$ here, since the only morphisms we need in the following subsections are \emph{auto}morphisms in $\sqcup \C$, for which we will now describe a simpler data structure:
Indeed, being an automorphism of $T$ forces morphisms between the transitive cofactors to remain in the same component.
This translates to the isomorphism
\begin{align} \label{eq:Aut}
  \Aut_{\sqcup \C}(T) \cong \Aut_\C(C_1) \wr S_{m_1} \times \cdots \times \Aut_\C(C_\ell) \wr S_{m_\ell} \mbox{,}
 \end{align}
 where each element is an $\ell$-tuple with $c$-th entry given by a pair
 \[
   ((\alpha_{c 1}, \ldots, \alpha_{c m_c}), \sigma_c) \in \Aut_\C(C_c) \wr S_{m_c} \coloneqq \Aut_\C(C_c)^{\times m_c} \rtimes S_{m_c} \mbox{.}
\]
In \Cref{sec:SchreierSim} we will describe the implementation of $\CoequalizerOfIdentityAndAutomorphisms$ in $\T \coloneqq \sqcup \C$ using the Schreier-Sims orbit algorithm.

\subsubsection{The category $\GFinSet$ as a tower.}

A \emph{skeletal} model of $\GFinSet$ is finally given by the coproduct completion $\GFinSet \simeq \sqcup \Coeq(\Gamma)$.

The automorphism group of an object $C_c \coloneqq \Gamma/U_c \in \Coeq(\Gamma)$ is $\Aut(\Gamma/U_c) \cong \Nor_\Gamma(U_c) / U_c$ by \eqref{eq:mor_in_Coeq}, and any automorphism of $\coprod_{c=1}^\ell (\Gamma/U_c)^{\sqcup m_c} \in \GFinSet$ can be modeled according to \eqref{eq:Aut} by an element in the group $\times_{c=1}^\ell \Nor_H(U_c) \wr S_{m_c}$, which is again an $\ell$-tuple with $c$-th entry given by a pair $((\nu_{c 1}, \ldots, \nu_{c m_c}), \sigma_c) \in \Nor_\Gamma(U_c) \wr S_{m_c} \leq \Gamma \wr S_{m_c} \coloneqq \Gamma^{\times m_c} \rtimes S_{m_c}$.

The finite products in $\sqcup \Coeq(\Gamma)$ can be easily computed using the table of Burnside marks $\{ \beta_{\Gamma / U_c} \}_{c=1}^\ell$ of $\Gamma$ with $\beta_{\Gamma / U_c} \coloneqq \Hom_{\Coeq(\Gamma)}(-, \Gamma / U_c)$.
In fact, $\GFinSet \simeq \sqcup \Coeq(\Gamma)$ is an elementary topos.
We have implemented this tower in the \CAP-based package \FinGSetsForCAP\footnote{The package contains a primitive implementation of the category $\FinSet\textnormal{-}\Gamma$ of right $\Gamma$-sets by the third author. We will replace this primitive implementation by a compilation using \CompilerForCAP \cite{CompilerForCAP,ZickgrafDoktor} of the tower $\sqcup \Coeq(\Gamma^\mathrm{op}) \simeq \FinSet\textnormal{-}\Gamma$.} \cite{FinGSetsForCAP} together with its universal property as a finite colimit completion.
This universal property is given by the composition of the ur-algorithm of the coequalizer completion (see \Cref{sec:SchreierSim}) and the ur-algorithm for the finite coproduct completion, the latter being formal.
The computational access to this universal property is key for the computation of the Kleisli composition in the biset category, as we have seen in \Cref{sec:Biset}.

\subsection{A categorical interpretation of the Schreier-Sims algorithm} \label{sec:SchreierSim}

In order to implement the Kleisli composition in \Cref{fig:Kleisli} we need to compute the special coequalizer $\CoequalizerOfIdentityAndAutomorphisms$ \eqref{eq:coequalizer_of_identity_and_automorphisms} in the category $\T = \sqcup \C$, where $\C$ is a skeletal object-finite category $\{C_1, \ldots, C_\ell\}$ equipped with such coequalizers, specialized to $\C = \Coeq(\Gamma)$ for $\Gamma = K$.

It turns out that computing $\coeq_T(\alpha)$ as in \eqref{eq:coequalizer_of_identity_and_automorphisms} for $T \in \T = \sqcup \C$ and $\alpha = (\alpha_1, \ldots, \alpha_k) \in \Aut_\T(T)^{\times k}$ reduces to the \emph{Schreier--Sims orbit algorithm} and gives the latter an elegant categorical interpretation:
Because all morphisms in the above tuple $\alpha$ are automorphisms we can reduce the computation of the coequalizer to each cofactor of $T \coloneqq \coprod_{c=1}^\ell C_c^{\sqcup m_c}$, and then use the functoriality of the coproduct to recover the coequalizer.
So we can assume $T$ to be a coproduct $T = C^{\sqcup m}$ of the same object $C \in \C$.
For $r=1, \ldots, k$ each $\alpha_r = ((\nu_{r 1}, \ldots, \nu_{r m}), \sigma_r) \in \Aut_\T(T) = \Aut_\C(C) \wr S_m$.
Define $P \coloneqq \langle \sigma_1, \ldots, \sigma_k \rangle \leq S_m$.
Without loss of generality we can assume that $P$ acts transitively on the set $\{1, \ldots, m\} \cong_P S \setminus\! P$ (from the right) with $S = \stab_P(1)$, otherwise we restrict to $P$-orbits.
Under these assumptions $D \coloneqq \coeq_T(\alpha) \in \C \subseteq \T$ is a factor object of $C$ in $\C$.
In order to compute $D \in \C$ and the $m$ projections $\pi_i\colon C \to D$ we invoke the Schreier--Sims orbit algorithm for the transitive permutation group $P$ and obtain $m$ words $t_1 = (), t_2, \ldots, t_m$ and $z \coloneqq m k - m + 1$ words $s_j$ in $k$ indeterminates such that the evaluation $t_i(\sigma_1, \ldots, \sigma_k) \in P$ is the transversal of the right coset $S t_i(\sigma_1, \ldots, \sigma_k) \subseteq P$ corresponding to the point $i \in \{1, \ldots, m\}$ and the evaluations $s_j(\sigma_1, \ldots, \sigma_k) \in P$ for $j=1, \ldots, z$ generate the stabilizer $S$.
We now replace the $\sigma_r$'s with the $\alpha_r$'s in the evaluations and applying them to $1 \in \{1, \ldots, m\}$.
Doing this, we first obtain $z$ many automorphisms $\delta_j \coloneqq 1.s_j(\alpha_1, \ldots, \alpha_k) \in \Aut_\C(C)$ with $D = \coeq_C( \delta )$ along with the coequalizer projection $\pi: C \to D$.
Finally, we obtain the $m$ projections as $\pi_i \coloneqq 1.t_i^{-1}(\alpha_1, \ldots, \alpha_k) \pi\colon C \to D$.

\subsection{Action pair as a data structure for a biset} \label{sec:action_pair}

According to \eqref{eq:biset} a morphism $X \in \FinBiset(G,H)$ corresponding to the $(H,G)$-biset ${}_H X_G$ is a functor $X\colon G \to \HFinSet \simeq \sqcup \Coeq(H)$.
We now describe the data structure for $X$.
Let $G$ be generated by $\{ g_a \mid a=1,\ldots, n \}$:
\[
  G = \langle g_a \mid a=1,\ldots, n \rangle \mbox{.}
\]
The data structure for $X$ is given by an action pair $\chi = \left({}_H X, ( \chi_a )_{a=1}^n \right)$ which simply encodes the corresponding group homomorphism $\chi\colon G \to \Aut({}_H X)$ from \eqref{eq:chi}, i.e., with $\chi_a \coloneqq \chi(g_a) \in \Aut({}_H X)$.

According to \Cref{sec:Gamma-FinSet_as_a_tower_implementation}, the first component of the pair $\chi$ is the left $H$-left ${}_H X \in \HFinSet$, which was represented above by an $\ell$-tuple $(m_1, \ldots, m_\ell)$ of multiplicities.
Furthermore, each $H$-equivariant automorphism $\chi_a$ of ${}_H X$ in the second component is given by an $\ell$-tuple
\[
  (\chi_{a, c})_{c=1}^\ell \coloneqq \left( ((\nu_{a, c 1}, \ldots, \nu_{a, c m_c}), \sigma_{a, c}) \right)_{c=1}^\ell
  \in
  \bigtimes_{c=1}^\ell \Nor_H(U_c) \wr S_{m_c} \mbox{,}
\]
where $\{ U_c \}_{c = 1}^\ell$ are representatives of the conjugacy classes of subgroups of $H$.
Note that the action of $G$ described in \eqref{Hxg} corresponds to the action of the subgroup $\times_{c=1}^\ell S_{m_c}$ generated by $(\sigma_{a,1}, \ldots, \sigma_{a, \ell})$ for $a = 1, \ldots, n$.

The additive decomposition of an action pair into an action pair of transitive bisets starts by decomposing an action pair $\chi = \left((m_c)_{c=1}^\ell, \left((\chi_{a, c})_{c=1}^\ell\right)_{a=1}^n\right)$ into its $\ell$ many components.
The $c$-th component can then be decomposed further by considering the orbits of the permutation group $P_c \coloneqq \langle \sigma_{a, c} \mid a = 1, \ldots, n \rangle_{S_{m_c}}$ on $\{1, \ldots, m_c\}$.


\subsubsection{From a subgroup to a transitive action pair.} \label{sec: subgroup to action pair}

As described in \Cref{sec:goursat_bouc} a subgroup $D \leq H \times G$ defines a transitive $(H,G)$-biset $X \coloneqq (H \times G) / D$, which we view as a morphism in $\mathbf{Biset}(G,H)$.
We now describe how to construct the corresponding action pair $\chi$:
To recover the multiplicities in the first entry of the action pair we decompose ${}_H X$ into $H$-orbits:
\[
  (H \times G) / D \cong_H \coprod_{i = 1}^m H \cdot (1_H, t_i)D \mbox{,}
\]
where $(t_1, \ldots, t_m )$ is a set of transversals of left cosets $G/P_2(D)$.
An easy calculation shows that $H \cdot (1_H, t_i)D \cong_H H / K_1(D)$ as left $H$-sets, for $i \in \{ 1, \ldots, m \}$.
We then find the number $c$ such that $H/K_1(D) \cong_H H/U_c$ and define the first entry of $\chi$ as the list of multiplicities $(m_1, \ldots, m_\ell)$ with $m_d = m \delta_{c, d}$.
It remains to describe the $c$-th component $\chi_{a, c} = ((\nu_{a, c 1}, \ldots, \nu_{a, c m_c}), \sigma_{a, c})$ of each $\ell$-tuple in the second entry for $a = 1, \ldots, n$.
For a transversal $t_i$ and and a generator $g_a$
\begin{align*}
  (1_H, t_i)D \cdot g_a &  = (1_H, g_a^{-1}t_i) D = ( 1_H, t_{\sigma_{a, c}(i)} g_{a, i} ) D \\
  & = ( h_{a, i}^{-1}, t_{\sigma_{a, c}(i) } ) D = h_{a, i}^{-1} \cdot (1_H, t_{\sigma_{a, c}(i)})D \mbox{,}
\end{align*}
where $t_{\sigma_{a, c}(i)}$ is the representative of the left coset containing $g_a^{-1} t_i$, $g_{a, i} \coloneqq ( t_{\sigma_{a, c}(i)} )^{-1} g_a^{-1} t_i \in P_2(D)$, and $h_{a, i} \in P_1(D)$ such that $(h_{a, i},g_{a, i}) \in D$.
The element $h_{a, i}^{-1}$ defines an $H$-equivariant iso $H \cdot (1_H, t_i) D \xrightarrow{\sim} H \cdot (1_H,t_{\sigma_{a, c}(i)})D$, which is independent of the choice of $h_{a, i}$.
We then obtain using \eqref{eq:mor_in_Coeq} an automorphism $\nu_{a, c i} \coloneqq h h_{a, i}^{-1} h^{-1}\colon H/U_c \to H/U_c$, where $h \in H$ is such that $U_c = {}^h K_1(D)$.

\subsubsection{From a transitive action pair to a subgroup.}

Let $\chi$ be the action pair of a transitive $(H,G)$-biset ${}_H X_G$.
Then ${}_H X_G$ is isomorphic to $(H\times G)/\stab_{H\times G}(x)$ for any element $x \in X$.
In order to compute the stabilizer $D = \stab_{H\times G}(x)$ we use the following fact, which is a reformulation of \eqref{eq:D}:
\begin{equation} \label{eq:stab_decomposition}
  \stab_{H \times G}(x) = \left\langle \ \stab_H(x)\times \{1_G\} \ \cup \ \{(\phi(g),g) \mid g \in \stab_G(H \cdot x)\} \ \right\rangle \mbox{,}
\end{equation}
with $\phi\colon \stab_G(H \cdot x) \to H$ being any map satisfying $\phi(g) \cdot x = x \cdot g$, where $\stab_G(H \cdot x)$ refers to the action \eqref{Hxg} of $G$ on the $H$-orbits of ${}_H X$.

We now compute $D = \stab_{H \times G}(x)$ for a particular choice of $x$:
Since the biset is transitive, the first entry of $\chi$ is the $\ell$-tuple $(0, \ldots, 0, m_c, 0, \ldots 0)$.
We choose $x$ to be the left coset $1_H U_c = U_c$ in the first cofactor $H / U_c$ of ${}_H X$. Hence $\stab_H(x) = U_c$.
Then we compute the subgroup $\stab_G( H \cdot x ) \leq G$ as the stabilizer of $1$ under the action of $G$ on $\{1, \ldots, m_c\}$ defined by $g_a \mapsto \sigma_{a, c} \in S_{m_c}$.
Now define $\phi\colon \stab_G( H \cdot x ) \to H$ on $g$ by the evaluation $\phi(g) \coloneqq w(\chi_{1,c}, \ldots, \chi_{n, c})(1)$, where $w_g$ is an abstract word in $n$ indeterminates such that $w_g(g_1, \ldots, g_n) = g$.


\section{The relation to the category of generalized morphism in the Abelian group case}

One can further restrict $\FinBiset$ to the full subcategory on finite \emph{Abelian} groups.
Mackey's formula \eqref{eq:Mackey} states that the composition of two transitive bisets of Abelian groups is at most a multiple of a transitive biset.
When ignoring these multiplicities in the composition (e.g., by setting them to $1$), the transitive bisets of Abelian groups define a category, which coincides with the category  $\mathbf{G}(\mathbf{FinAb})$ of generalized morphisms of the category $\mathbf{FinAb}$ of finite Abelian groups.
The category of generalized morphisms $\mathbf{G}(\mathbf{A})$ of a general Abelian category $\mathbf{A}$ was introduced in \cite{BL_GabrielMorphisms} and further investigated and implemented in \CAP in \cite{GutscheDoktor}.

While $\FinBiset$ is enriched over the category of finitely generated free commutative monoids (where the pre-additive structure is a manifestation of the coproduct of $\Gamma$-sets), $\mathbf{G}(\mathbf{A})$ is enriched over commutative inverse monoids \cite{BL_GabrielMorphisms}.

\begin{credits}
\subsubsection{\ackname}
This work is a contribution to Project `Compute the equivariant Orlik-Solomon algebra of a matroid' of SPP 2458 ``Combinatorial Synergies'', funded by the Deutsche Forschungsgemeinschaft (DFG, German Research Foundation).
The third author was supported by the ``Qompiler Project'', funded by the German Federal Ministry for Economic Affairs and Climate Action to develop \CompilerForCAP \cite{CompilerForCAP,ZickgrafDoktor}, which is used to compile the implementation described in this paper.
The second author thanks his PhD supervisor Dr.~Radu Stancu for extensive discussions and continuous support.
We are grateful to both anonymous reviewers for their valuable suggestions.
\end{credits}

\bibliographystyle{splncs04}
\bibliography{Bisets}

@misc{web2023,
      title={Biset functors for categories}, 
      author={Peter Webb},
      year={2023},
      eprint={2304.06863},
      archivePrefix={arXiv},
      primaryClass={math.RT},
      note={(\href{https://arxiv.org/abs/2304.06863}{arXiv:2304.06863})},
}

@inproceedings{GPSSyntax,
    author = {Gutsche, Sebastian and Posur, Sebastian and Skarts{\ae}terhagen, {\O}ystein},
     title = {\href{https://ceur-ws.org/Vol-2307/paper21.pdf}{On the Syntax and Semantics of $\mathtt{CAP}$}},
 booktitle = {O. Hasan, M. Pfeiffer, G. D. Reis (eds.): Proceedings of the Workshop Computer Algebra in the Age of Types, Hagenberg, Austria, 17-Aug-2018},
      year = {2018},
}

@manual{FinGSetsForCAP,
author = {Barakat, Mohamed and Mickisch, Julia and Talleux, Marc and Zickgraf, Fabian},
title = {{$\mathtt{FinGSetsForCAP}$ -- The (skeletal) elementary topos of finite G-sets}},
year = {2017--2026},
url = {https://homalg-project.github.io/pkg/FinGSetsForCAP},
}

@manual{FiniteCocompletions,
author = {Barakat, Mohamed and Talleux, Marc},
title = {{$\mathtt{FiniteCocompletions}$ -- Finite (co)product/(co)limit (co)completions}},
year = {2022--2026},
url = {https://homalg-project.github.io/pkg/FiniteCocompletions},
}

@manual{CompilerForCAP,
author = {Zickgraf, Fabian},
title = {{$\mathtt{CompilerForCAP}$ -- Speed up and verify categorical algorithms}},
year = {2020--2026},
url = {https://homalg-project.github.io/pkg/CompilerForCAP},
}

@phdthesis{ZickgrafDoktor,
    author = {Zickgraf, Fabian},
     title = {\href{https://dx.doi.org/10.25819/ubsi/10541}{CompilerForCAP -- Building and compiling categorical towers in algorithmic category theory}},
    school = {University of Siegen},
      year = {2024},
      type = {Dissertation},
}

@book{bouc2010biset,
    AUTHOR = {Bouc, Serge},
     TITLE = {\href{https://doi.org/10.1007/978-3-642-11297-3}{Biset functors for finite groups}},
    SERIES = {Lecture Notes in Mathematics},
    VOLUME = {1990},
 PUBLISHER = {Springer-Verlag, Berlin},
      YEAR = {2010},
     PAGES = {x+299},
      ISBN = {978-3-642-11296-6},
   MRCLASS = {20J15 (19A22 20C15 20C20)},
  MRNUMBER = {2598185},
MRREVIEWER = {Robert\ Hartmann},
}

@unpublished{BL_GabrielMorphisms,
    AUTHOR = {Barakat, Mohamed and Lange-Hegermann, Markus},
     TITLE = {{G}abriel morphisms and the computability of {S}erre quotients with applications to coherent sheaves},
      year = {2014},
      note = {(\href{https://arxiv.org/abs/1409.2028}{arXiv:1409.2028})},
}

@phdthesis{GutscheDoktor,
    author = {Gutsche, Sebastian},
     title = {\href{https://nbn-resolving.org/urn:nbn:de:hbz:467-12411}{Constructive category theory and applications to algebraic geometry}},
    school = {University of Siegen},
      year = {2017},
      type = {Dissertation},
}

@Misc{GP_Rundbrief,
   author = {Gutsche, Sebastian and Posur, Sebastian},
    title = {\href{https://fachgruppe-computeralgebra.de/data/CA-Rundbrief/car64.pdf}{CAP: categories, algorithms, programming}},
    PAGES = {14--17},
     year = {2019},
   howpublished ={{\em Computeralgebra-Rundbrief}, \textbf{64}, 14--17, March},
}

@incollection {PosCCT+arXiv,
    AUTHOR = {Posur, Sebastian},
     TITLE = {\href{https://arxiv.org/abs/1908.04132}{Methods of constructive category theory}},
 BOOKTITLE = {Representations of {A}lgebras, {G}eometry and {P}hysics},
    SERIES = {Contemp. Math.},
    VOLUME = {769},
     PAGES = {157--208},
 PUBLISHER = {Amer. Math. Soc., [Providence], RI},
      YEAR = {2021},
   MRCLASS = {18E10 (18A25 18E05 18E25)},
  MRNUMBER = {4254099},
}

\end{document}